\documentclass[12, eqno]{amsart}
\usepackage{amsmath}
\usepackage{amscd,amsthm,amsfonts,amsopn,amssymb,verbatim}
\parskip =\medskipamount
\numberwithin{equation}{section}
\newtheorem{theorem}{Theorem} 

\theoremstyle{remark}
\DeclareMathOperator{\supp}{supp\,}

\def\dist{\rm dist}
\def\be{\begin{equation}}
\def\ee{\end{equation}}
\def\vp{\varphi}
\def\ve{\varepsilon}

\allowdisplaybreaks

\begin{document}

\title{On pair correlation for generic diagonal forms}
\date{\today}
\author{J.~Bourgain}
\address{J.~Bourgain, Institute for Advanced Study, Princeton, NJ 08540}
\email{bourgain@math.ias.edu}
\thanks{The author was partially supported by NSF grants DMS-1301619}

\begin{abstract}
We establish new pair correlation results for certain generic homogenous diagonal forms evaluated on the integers. 
Methods are analytic leading to explicit quantitative statements.
\end{abstract}
\maketitle
\section {Introduction and statements}

In \cite {S}, an analytic and quantitative approach to the pair correlation problem for (measure) generic binary 
quadratic forms $\alpha m^2+mn+ \beta n^2$, $\alpha, \beta>0)$ is given, establishing Poisson behavior.
The method used in \cite {S}  has been developed further by several authors, in particular in \cite {V}, to which we will refer later.
On the other hand, Sarnak's argument as it stands does not seem to apply to diagonal forms $m^2+\beta n^2, \beta>0$, corresponding to the Laplace eigenvalues of a
rectangular billard, due to lack of parameters.
Pair correlation results for binary quadratic forms have been obtained in \cite{EMM}, based on ergodic methods \big(in this case the pair correlation statistics amounts to the
distribution of quadratic forms of (2, 2)-signature\big).
A considerable advantage of the ergodic method is to provide deterministic results, with the quadratic forms being subject to a diophantine assumption.
Those results are qualitative and only weak quantitative statements seem extractable from this technology.
The purpose of this Note is to provide an alternative for Sarnak's approach which applies in the diagonal case and also in situations where no dynamical treatment is
known.
The technique is closely related to arguments in \cite{BBRR} and \cite{B}.

A first model, suggested to the author by Z.~Rudnick, is that of a generic positive definite diagonal ternary quadratic form
$$
Q(\bar x) =Q (x_1, x_2, x_3)= x_1^2 +\alpha x_2^2+\beta x_3^2, \alpha, \beta >0
$$
Restricting the variables to $\mathbb Z\cap [0, N]$, typical spacings are expected to be $O(\frac 1N)$ so that the natural interval size in the pair correlation  problem
for (1.1) is $O(\frac 1N)$.
For $T$ large and $a<b$, define
$$
\begin{aligned}
R(a, b; T)= & \frac 1{T^{3/2}}|\{(\bar m, \bar n)\in \mathbb Z_+^3\times \mathbb Z_+^3; \bar m\not= \bar n, Q(\bar m)\leq T, Q(\bar n)\leq T \text { and } \\[8pt]
&Q(\bar m)-Q(\bar n)\in [a, b]\}|
\end{aligned}
\eqno{(1.1)}
$$
and
$$
c=\lim_{\ve\to 0}\frac 1\ve \Big|\Big\{(x, y)\in\mathbb R^3_+\times\mathbb R_+^3; Q(x)\leq 1, |Q(x)-Q(y)|<\frac\ve 2\Big\}\Big|.\eqno{(1.2)}
$$
\begin{theorem}
Let $Q=Q_{\alpha, \beta}$ be defined by (1.1) with $\alpha,\beta \in [\frac 12, 1]$ parameters.

Fix $\frac 12<\rho<1$.
Then almost surely in $\alpha, \beta$ the following statement holds.
Let $T\to\infty$ and $a<b$, $|a|, |b|<O(1)$ and $T^{-\rho}< b-a<1$.
Then
$$
R(a, b; T)\sim c T^{\frac 12}(b-a).\eqno{(1.3)}
$$
\end{theorem}

The second result relates to Vanderkam's work \cite {V} on the pair correlation for homogenous degree $k$ forms in $k$ variables (and which is an extension of \cite {S}).
While in \cite {V} a measure generic result for the full space of such forms is established, we restrict ourselves to diagonal forms, proving a similar result.
Define for given $k\geq 3$
$$
F(x_1, \ldots, x_k)= x_1^k+\alpha_2x_2^k+\cdots+ \alpha_kx_k^k, \text { where} \ \alpha_2, \ldots, \alpha_k\in \Big[\frac 12, 1\Big]\eqno{(1.4)}
$$
$$
c=|\{x\in\mathbb R_+^k; F(x)\leq 1\}|\eqno{(1.5)}
$$
$$
R(a, b, T)=\frac 1T|\{(\bar m, \bar n)\in \mathbb Z_+^k\times \mathbb Z^k_+; \bar m\not=\bar n, F(\bar m)\leq T, F(\bar n)\leq T \text { and }
F(\bar m)-F(\bar n)\in [a. b]\}|.
$$

\begin{theorem}
Let $F$ be defined by (1.4) with $\alpha_2, \ldots, \alpha_k\in [\frac 12, 1]$ parameters.
There is some $\rho>0$ such that almost surely in $\alpha_2, \ldots\alpha_k$ the following statement holds.
Let $T\to\infty$ and $0<a< b<O(1)$ and $b-a> T^{-\rho}$. Then
$$
R(a, b, T)\sim c^2(b-a)\eqno{(1.6)}
$$
\end{theorem}

More precise statements will follow from the arguments below.
Some steps were not made quantitatively explicit in order not to over-complicate the exposition.
It turns out that in the diagonal case, the most direct harmonic analysis attack introducing Gauss sums does not seem to succeed.
Instead, we follow a slightly different approach based on distributional properties of certain Dirichlet sums (cf. \cite{BBRR} and \cite {B}).

\section
{Preliminaries}

We make a few comments/reductions related to Theorem 1 (a similar discussion holds for Theorem 2).
Clearly, Theorem 1 is equivalent to the statement (by rescaling)
$$
\frac 1{T^{\frac 12}} R(a, b; T)\sim\frac {b-a}{2.T^{8/3}} |\{(\bar x, \bar y)\in\mathbb R^3_+\times \mathbb R^3_+; Q(\bar x)\leq T, |Q(\bar x)-Q(\bar y)|< T^{2/3}\}|.
\eqno{(2.1)}
$$

We make a localization of the variables $x_i \in I_i =[u_i -\Delta u_i, u_i+\Delta u_i]$ and $y_i \in J_i =[v_i-\Delta v_i, v_i +\Delta v_i]$
(where $I_i, J_i\subset \mathbb R_+$, hence $I_i, J_i\subset [0, cT^{\frac 12}])$ and such that $u_i> T^{\frac 12-\ve}$, $\Delta u_i=T^{\frac 12-3\ve}$ and similarly
for $v_i$ and $\Delta v_i$.
Then (2.1) will follow from

$$
\begin{aligned}
&|\{ (\bar m, \bar n)\in\mathbb Z^3\times\mathbb Z^3; \bar m\in I_1\times I_2\times I_3, \bar n \in J_1\times J_2\times J_3, \bar m \not= \bar n \text { and } Q(\bar m)- Q(\bar n)\in [a,
b]\}|\sim\\[10pt]
& \frac {b-a}{2T^{2/3}} |\{ (\bar x, \bar y) \in \mathbb R^3\times\mathbb R^3; \bar x\in I_1\times I_2\times I_3, \bar y \in J_1\times J_2\times J_3 \text { and } |Q(\bar x)-
Q(\bar y)|< T^{2/3}\}|
\end{aligned}
\eqno{(2.2)}
$$
We may further ensure that (as we explain next)
$$
\dist(I_i, J_i) \sim |u_i-v_i| >T^{\frac 12 -\ve}.\eqno{(2.3)}
$$
Since $Q (\bar m) -Q(\bar n) = m^2_1 -n_1^2+\alpha(m^2_2 -n_2^2)+ 
\beta (m_3^2 -n_3^2)$ and $b-a<1$, the contribution of say $|m_2-n_2|<T ^{\frac 12+\ve}$ to $R(a, b; T)$ is at most (denoting $\xi =\frac
{a+b}2$).

$$
T^{-3/2} \Bigg| \Bigg\{ 
\begin{aligned}
(m_2, n_2, m_3, n_3)\in \mathbb Z^4_+ ; m_i,  n_i\lesssim & T^{\frac 12}, m_3\not= n_3, |m_2-n_2|<T^{\frac 12-\ve} \text {and}\\[10pt]
& \Vert \alpha (m_2^2 -n_2^2) +\beta (m_3^2 -n_3^2)- \xi\Vert < b-a\end{aligned}\Bigg\}
\Bigg|\eqno{(2.4)}
$$
$$
+T^{-3/2} \Big|\Big\{ (\bar m, \bar n)\in \mathbb Z_+^6; m_i, n_i \lesssim T^{\frac 12}, m_3 =n_3, \bar m \not =\bar n \text { and } |m^2_1-n_1^2 +\alpha (m_2^2-n_2^2)
-\xi|< b-a\Big\}\Big|.
\eqno{(2.5)}
$$
Since $\alpha, \beta$ are measure generic, hence diophantine, (2.4) can be estimated by
$$
T^{-3/2} T^{1-\ve} \big(T^{1+\frac \ve 2} (b-a) +1\big) <T^{\frac 12-\frac \ve 2} (b-a) + T^{-\frac 12-\ve} <T^{\frac 12-\frac \ve2}(b-a).
$$
For (2.5), since $|\xi|<O(1)$ and $\bar m \not= \bar n$, either $m_2\not= n_2$ or $m_2 =n_2, m_3= n_3, m_1, n_1< O(1)$.

Again the $m_2\not = n_2$ contribution is at most
$$
T^{-\frac 32} T^{\frac 12} T^{1+\frac\ve 2} (b-a) +1\big)< T^{\frac \ve 2}(b-a) + T^{-1} < T^{\frac \ve 2} (b-a).
$$
In order to evaluate the contribution of $|m_1-n_1|< T^{\frac 12-\ve}$, simply replace $Q(x)$ by $\frac 1\beta x_1^2+\frac \alpha \beta x_2^2 +x_3^2$ and exploit that 
$\frac 1\beta$, $\frac \alpha\beta$ are diophantine.
This justifies (2.3).

Next, observe that in order to establish Theorem 1, we need to consider the $T^{\rho+\ve}$ cases for the interval $[a, b]$.
This means that the measure of the exceptional parameter set in $(\alpha, \beta)$ for a fixed interval $[a, b]$ has to be multiplied by $T^{\rho+\ve}$.

Similar considerations apply to Theorem 2, replacing $T^{\frac 12}$ by $T^{\frac 1k}$ and variables range $0<x_i, y_i \lesssim T^{1/k}$.
Condition (2.3) is replaced by
$$
\dist (I_i, J_i)\sim |u_i-v_i|> T^{\frac 1k-\ve}\eqno{(2.6)}
$$
which we justify, as the argument differs.
Assume $|m_1-n_1|<T^{\frac 1k-\ve}$.

Note first that for a fixed interval $[a, b]$, we have for fixed $\bar m\not=\bar n$

$$
\begin{aligned}
&\Big|\Big\{ (\alpha_2, \ldots, \alpha_k) \in \Big[\frac 12, 1\Big]^{k-1}; |m_1^k- n_1^k +\alpha_2 (m_2^k-n_2^k)+\cdots+ \alpha_k (m_k^k -n_k^k)-\xi|
< b-a\Big\}\Big|\\[10pt]
&\lesssim \Big(\sum^k_{i=1} |m_i^k- n_i^k|\Big)^{-1}.
\end{aligned}
\eqno{(2.7)}
$$
Summing (2.7) over $\bar m, \bar n$ and using the geometric/arithmetical mean inequality, the contribution of $m_1\not= n_1$ to $R_1 (a, b; T)$ is at most

$$
\begin{aligned}
&\frac C T \sum_{\substack{ 1\leq |m_1-n_1|< T^{\frac 1k-\ve}\\ m_i, n_i\leq T^{\frac 1k} (2\leq i\leq k)}} \ \frac 1{\sum^k_{i=1} |m_i^k -n_i^k|}\leq\\[10pt]
&\frac CT \Big(\iint_{\substack{1\leq u\leq T^{\frac 1k-\ve}\\ v\leq T^{\frac 1k}}} \ \frac {dudv}{u^{\frac 1k}\,  v^{\frac {k-1}k}}\Big)
\Big( T^{\frac 1k}+\iint_{1\leq u\leq v\leq T^{\frac 1k}} \ \frac {dudv}{u^{\frac 1k} v^{\frac{k-1} k}}\Big)^{k-1} \leq\\[10pt]
&\frac CT T^{(\frac 1k-\ve)\frac {k-1}k} T^{\frac {1}{k^2}} T^{\frac{k-1}k} \lesssim T^{-\ve\frac {k-1}k}.
\end{aligned}
\eqno{(2.8)}
$$

If $m_1=n_1$, then say $m_2\not= n_2$ contributing

$$
\begin{aligned}
&C\frac {T^{\frac 1k}}T \sum_{\substack{m_i, n_i \leq T^{\frac 1k} (2\leq i\leq k)\\ m_2\not= n_2}} \ \frac 1{\sum^k_{i=2} |m_i^k -n_i^k|}\leq\\[10pt]
&C T^{\frac 1k-1} \Big( \iint_{1\leq u\leq v\leq T^{\frac 1k}} \ \frac {dudv}{u^{\frac 1{k-1}} v}\Big)
\Big( T^{\frac 1k}+\iint_{1\leq u\leq v\leq T^{\frac 1k}}
\ \frac {dudv}{u^{\frac 1{k-1}} v} \Big)^{k-2}\leq\\[10pt]
& CT^{\frac 1k-1} (\log T) T^{\frac {k-2}{k(k-1)}} T^{\frac {k-2}k} \lesssim T^{-\frac 1{k(k-1)}} \log T.\end{aligned}
\eqno{(2.9)}
$$
Multiplying (2.8), (2.9) with the number of intervals $[a, b]$ under consideration, i.e. $O(T^\rho)$, the argument is conclusive letting $\rho =\min \big(\frac \ve 4, \frac
1{2k(k-1)}\big) =\frac\ve 4$.
This justifies (2.6).

\section
{Proof of Theorem 1}

We establish (2.2) with the intervals $I_i, J_i(1\leq i\leq 3)$ introduced as above and in particular the separation condition (2.3).

Let $\xi =\frac {a+b}2 , \delta =\frac {b-a}2$.
Let $I_3=[u-\Delta u, u+\Delta u], I_3 =[v-\Delta v, v+\Delta v], k_0= u^2-v^2$ and note that by (2.3)
$$
|k_0|\geq |u-v|^2 > T^{1-2\ve}\eqno{(3.1)}
$$
while for $x\in I_3, y\in J_3$
$$
|x^2 -y^2-k_0|\lesssim T^{\frac 12} (\Delta u +\Delta v) < T^{1-3\ve}.\eqno{(3.2)}
$$
Rewrite $Q(\bar m) -Q(\bar n)\in [a, b]$ as
$$
|m_1^2-n_1^2 +\alpha (m_2^2- n_2^2)+\beta (m_3^2 -n_3^2)-\xi |<\delta.\eqno{(3.3)}
$$
Taking into account (3.1), (3.2) and  exploiting a usual upper/lower bounding argument, (3.3) may be replaced by
$$
\Big|\frac {m_1^2 - n_1^2 +\alpha (m_2^2 -n_2^2)-\xi}{m_3^2- n_3^2} +\beta\Big| < \frac \delta{|k_0|}\eqno{(3.4)}
$$
or, assuming $k_0 >0$, i.e. $u>v$ and taking into account the variable localization

$$
|\log (n_1^2 -m^2_1+ \alpha (n_2^2 - m^2_2)+\xi)-\log (m_3^2-n_3^2)-\log \beta|<\frac \delta{\beta k_0}.\eqno{(3.5)}
$$
At this point, we use Fourier analysis.

Denote
$$
\qquad S_1 (t)= \sum_{\substack{m_i\in I_i, n_i\in J_i\\ i=1, 2}} \big(n_1^2 -m_1^2 +\alpha(n_2^2-m_2^2)+\xi\big)^{it}\eqno{(3.6)}
$$
$$
S_2(t)=\sum_{m_3 \in I_3, n_3\in J_3} (m_3^2 -n_3^2)^{it}\qquad\qquad .\eqno{(3.7)}
$$
Using (3.5), the l.h.s. of (2.2) may be expressed as
$$
\int S_1(t)\, \overline{S_2(t)} \ e^{-it\log \beta} \, \widehat {1_{[-\frac 1B, \frac 1B]}} (t) dt \text { with } B=\frac{\beta{k_0}}\delta.\eqno{(3.8)}
$$
Define $B_0=\frac {\beta k_0} {T^{2/3}} > T^{\frac 16}$ and split (3.8) as
$$
\begin{aligned}
&\frac {B_0}B\int S_1(t) \,  \overline {S_2(t)} e^{it\log\beta} \, \widehat{ 1_{[-\frac 1{B_0}, \frac 1{B_0}]}} (t) dt +\int S_1 (t) \overline {S_2(t)} e^{-it \log B} \big[
\widehat { 1_{[-\frac 1 B, \frac 1B]}} (t) -\frac {B_0}B \widehat{ 1_{[-\frac 1{B_0}, \frac 1{B_0}]}} (t)\big] dt\\[8pt]
&= (3.9)+(3.10).
\end{aligned}
$$

Clearly (3.9) amounts to
$$
\begin{aligned}
&\frac {B_0}{B} \Big|\Big\{ (\bar m, \bar n)\in \mathbb Z^3 \times\mathbb Z^3; m_i\in I_i, n_i \in J_i \text { and } \Big|\log
\frac {n_1^2 -m_1^2 +\alpha (n_2^2 -m_2^2)+\xi}{\beta(m_3^2-n_3^2)}\Big| < \frac 1{B_0}\Big\}\Big|\\
&\sim \delta T^{-\frac 23} \Big|\Big\{ (\bar m , \bar n) \in \mathbb Z^3\times\mathbb Z^3; m_i\in I_i, n_i\in J_i \text { and }
|m_1^2-n^2_1+\alpha (m_2^2 -n_2^2)+\beta (m_3^2 -n_3^2)-\xi |< T^{\frac 23}\Big|\Big\}
\end{aligned}
$$
and since $\xi =O(1)$
$$
\sim \delta T^{-\frac 23} \Big|\Big\{ (\bar x, \bar y)\in \mathbb R^3\times\mathbb R^3; x_i\in I_i, y_i \in J_j \text { and } |Q(x) -Q(y)|<T^{\frac 32}\Big\}\Big|\eqno{(3.11)}
$$
which is (2.2).
Note that so far, no further restrictions on $\alpha, \beta$ were made.

Since clearly
$$
(3.11) >\delta T^{2-15\ve}\eqno{(3.12)}
$$
it remains to ensure that
$$
|(3.10) |< \delta T^{2-16 \ve}\eqno{(3.13)}
$$
which will involve an additional parameter restriction.

Analyzing the expression $\widehat{1_{[-\frac 1B, \frac 1B]}} (t) -\frac {B_0}{B} \widehat {1_{[-\frac 1{B_0}, \frac 1{B_0}]}} (t) $ and relying on $L^2$-theory, one verifies that
$$
\begin{aligned}
&\Vert(3.10)\Vert^2_{L^2_\beta}\lesssim\\[8pt]
&\Big(\frac \delta{k_0}\Big)^2 \Big\{\int_{[|t|<T^{-2/3} k_0]} (|t|T^{2/3}k_0^{-1})^4 |S_1(t)|^2 |S_2(t)|^2+
\int_{[|t|> T^{-2/3}k_0]} \min \Big(1, \frac {k_0}{|t|\delta}\Big)^2 |S_1(t)^2 |S_2(t)|^2\Big\} \\[8pt]
\end{aligned}\eqno{(3.14)}
$$
$$
\lesssim \Big(\frac \delta{k_0}\Big)^2\Big\{ T^{-\frac 43+C\ve} \int_{[|t|< T^\ve]} t^4 |S_1|^2 |S_2|^2 +\int_{[|t|>T^\ve]} 
\min \Big( 1, \frac T{|t|\delta}\Big)^2 |S_1|^2 |S_2|^2\Big\}.
\eqno{(3.15)}
$$
Recall that $S_1$ also depends on $\alpha$, so that the measure of the $(\alpha, \beta)$-set where (3.13) fails is bounded by
$$
T^{-6+C\ve} \Big\{ T^{-\frac 43} \int_{[|t|<T^\ve} t^4 (Av_\alpha |S_1 (t)|^2) |S_2(t)|^2 +\int_{|t|>T^\ve]} \min \Big( 1, 
\frac T{|t|S}\Big)^2 (Av_\alpha  |S_1(t)|^2) |S_2(t)|^2\Big\}.
\eqno{(3.16)}
$$
It remains to bound (3.16). Write
$$
|S_1(t)|^2 =\sum_{\substack{m_i, r_i \in I_i; n_i, s_i \in J_i\\ (i=1, 2)}}
e^{it[\log\big(n_1^2- m^2_1 +\alpha (n^2_2-m^2_2)+\xi\big)-\log \big(s_1^2-r_1^2+\alpha (s_2^2-r^2_2)+\xi)]}\eqno{(3.17)}
$$
and average in $\alpha$.  The $\alpha$-derivative of the phase function equals
$$
\frac {n_2^2 -m_2^2}{ n^2_1 -m_1^2+\alpha (n_2^2 -m_2^2)+\xi} -\frac {s_2^2 -r_2^2}{s_1^2 -r^2_1+\alpha (s_2^2-r_2^2)+\xi}
\geq \frac {(n_2^2 -m_2^2)(s_1^2-r^2_1+\xi)-(s_2^2 -r^2_2)(n_1^2-m_1^2+\xi)}{T^2}.
$$
This allows to bound the $\alpha$-average of (3.17) by

$$
\begin{aligned}
&\frac {T^2}{|t|} \ \sum_{\substack{ m_i, r_i \in I_i; n_i, s_i\in J_i\\ i=1, 2}}
[1+|(n^2_2-m_2^2)(s_1^2-r_1^2+\xi) -(s_2^2-r_2^2)(n_1^2-m_1^2+\xi)|]^{-1} \leq\\[10pt]
&\frac {T^2}{|t|} \ \sum_{\substack{z_i, w_i\in\mathbb Z\\ T^{1-2\ve} <|z_i|, w_i|<T}} [1+|z_2 (w_1+\xi) - w_2(z_1+\xi)|]^{-1}\leq\\[10pt]
&\frac {T^2}{|t|} \log T \max_{k\in\mathbb Z} \big|\{ (z_1, z_2, w_1, w_2) \in\mathbb Z^4; T^{1-2\ve}< |z_i|, |w_i| <T \text { and }
z_2 (w_1+\xi)-w_2 (z_1+\xi)= k+O(1)\big\}\big|\\[10pt]
&\lesssim \log T.\frac {T^2}{|t|} \sum_{T^{1-2\ve} < z_2, w_2<T} \Big(1+\frac T{z_2} (z_2, w_2)\Big) < \frac {T^{4+\ve}}{|t|}.
\end{aligned}\eqno{(3.18)}
$$
By (3.18),

$$
(3.16) <T^{-4/3+C\ve} +T^{-2+C\ve} \max_{2^k>T^\ve} \Big\{ 2^{-k} \min \Big(1, 2^{-k}\frac T\delta\Big) \int_{[|t|\sim 2^k]}|S_2(t)|^2\Big\}\eqno{(3.19)}
$$
where for the first term we used the trivial bound on $S_2$.

The main point is the argument in the treatment of $S_2$.
Recall the definition (3.7).
Instead of restricting $m_3\in I_3, n_3\in J_3$, it will be convenient to use a smoother localization, replacing $1_{[u-\Delta u, u+\Delta u]}$ by $\vp \big(\frac {x-u}{\Delta u}\big)$
with $\vp$ a bumpfunction satisfying $0\leq \vp\leq 1$, $\supp \vp\subset [-1, 1], \vp=1$ on $[-1+T^{-\ve}, 1-T^{-\ve}]$.
Then, writing $m_3^2-n_3^2= (m_3-n_3)(m_3+n_3)=mn$
$$
S_2(t) =\sum_{m, n} \vp\Big(\frac {m+n-2u}{2\Delta u}\big)\vp \Big( \frac {m-n-2v}{2\Delta v}\Big) m^{it} n^{it}
$$
and the latter may be bounded (up to some factor $T^{C\ve}$) by expressions of the from
$$
\Big(\sum_{m\in M} m^{it}\Big)\Big(\sum_{n\in N} n^{it}\Big)
$$
with $M, N$ subintervals of $[1, T^{\frac 12}]$ (recall that $m_3-n_3 > T^{\frac 12-\ve})$.

A well-known argument involving the Mellin transform (see for instance \cite{BBRR}) permits then to replace essentially $\frac 1{\sqrt K}\Big|\sum_{1\leq k\leq K} k^{it}|$ by the Riemann
zeta function $\zeta \big(\frac 12+it')$ with $t'$ a translate of $t$ by $O\big((\log T)^2\big)$.
Consequently

$$
\begin{aligned}
2^{-\ell}\int_{2^\ell <t< 2^{\ell+1}} |S_2 (t)|^2 dt & < T^{C\ve} 2^{-\ell} \max_{K<T^{\frac 12}}\int_{2^\ell <t<2^{\ell+1}} \Big|\sum _{1\leq k\leq K}
k^{it}\Big|^4 dt\\[8pt]
&<T^{C\ve} 2^{-\ell} T\int_{[1<t<2^\ell + (\log T)^2]}\Big|\zeta \Big(\frac 12+it\Big)\Big|^4 < T^{1+C\ve} 2^{\ve\ell}
\end{aligned}\eqno{(3.20)}
$$
using the classical bound

$$
\int_1^{t_*}\Big|\zeta \Big(\frac 12+ it\Big)\Big|^4 dt\ll t^{1+\ve}_* \ \text { for all } \ t_* >1, \ve>0.\eqno{(3.21)}
$$
Finally, substituting (3.20) in (3.19) gives a measure estimate in the $(\alpha, \beta)$ parameter set of the form $T^{-1+C\ve}$ (where $\ve>0$ may be taken arbitrarily small).
As pointed out earlier, this measure needs to be multiplied with the number $T^{\rho+\ve}$ of intervals $[a, b]$ under consideration and Theorem 1 follows.

\section
{Proof of Theorem 2}

We follow a similar procedure as for Theorem 1, though we need less optimal estimates due to the fact that the number of intervals under consideration is only $O(1)$.

Recall the definition of the intervals $I_i, J_i$ $(1\leq i\leq k)$ and the separation condition (2.6).
It will suffice to show that for some $0<\kappa<1$ (depending on $k$)

$$
\begin{aligned}
&\big|\big\{ (\bar m, \bar n) \in\mathbb Z^k\times\mathbb Z^k; m_i\in I_i, n_i\in J_i \ \text { and } \ F(\bar m)- F(\bar n)\in [a, b]\big\}\big|\sim\\[8pt]
& \frac {b-a}2 \, T^{-1+\kappa}\big|\big\{(\bar x, \bar y)\in \mathbb R^k \times\mathbb R^k; x_i\in I_i, y_i \in J_i \text { and } |F(x) -F(y)|< T^{1-\kappa}\big\}\big|.
\end{aligned}
\eqno{(4.1)}
$$
Define

$$
S_1(t) =\sum_{\substack{ n_i\in I_i, n_i\in J_i\\ i=1, \ldots, k-1}} (n^k_1-m_1^k +\alpha_2(n_2^k -m_2^k)+\cdots+ \alpha_{k-1}
(n^k_{k-1} - m^k_{k-1}) \big)^{it}\eqno{(4.2)}
$$
$$
S_2(t)=\sum_{m\in I_k, n\in J_k} (m^k -n^k)^{it}\eqno{(4.3)}
$$
and evaluate

$$
\int S_1(t) \overline {S_2(t)}\, e^{-it \log \alpha_k} \widehat {1_{[-\frac 1B, \frac 1B]}} (t) dt\eqno{(4.4)}
$$
where $B= \frac {\alpha_k k_o}\delta, k_0\sim m_k^k -n_k^k$.

Let $B_0 =\frac {\alpha_k k_0}{T^{1-\kappa}} > T^{\kappa/2}$ and decompose (4.4) as

$$
\begin{aligned}
&\frac {B_0}{B} \int S_1(t)\, \overline{S_2(t)} e^{-it\log \alpha_k} \widehat {1_{[-\frac 1{B_0}, \frac 1{B_0}]}} (t) dt+\int S_1(t) \, \overline {S_2(t)}\,
e^{-it\log \alpha_k} \big[\widehat {1_{[-\frac 1B, \frac 1B]}} (t) -\frac {B_0}{B} \widehat {1_{[-\frac 1{B_0}, \frac 1{B_0}]}} (t)\big] dt\\[8pt]
& = (4.5)+(4.6)
\end{aligned}
$$
Then, taking $\kappa<\frac 1k$, (4.5) amounts again to

$$
\delta T^{-1+\kappa}\big|\big\{ (\bar x, \bar y)\in \mathbb R^k \times \mathbb R^k; x_i\in I_i, y_i\in J_i \ \text { and } \ |F(\bar x) -F(\bar y)|<
T^{1-\kappa}\big\}\big|
$$
which is (4.1) and at least $T^{1-C\ve}$ (recall that $\delta =O(1)$ in this case).

Thus we need to ensure that
$$
(4.6)< T^{1-C'\ve}\eqno{(4.7)}
$$
achieved by additional parameter restriction.

Using the $L^2_{\alpha_k}$-norm as before,

$$
\begin{aligned}
\Vert(4.6)\Vert^2_{L^2_{\alpha_k}} &\lesssim \Big(\frac \delta{k_0}\Big)^2 \Big\{\int_{[|t|<\frac {k_0}{T^{1-\kappa}}]}
\Big(\frac {|t|T^{1-\kappa}}{k_0}\Big)^4 |S_1|^2 \, |S_2|^2 +\int_{[|t|>\frac {k_0}{T^{1-\kappa}}]} \min \Big(1, \frac {k_0}{\delta t}\Big)^2
|S_1|^2 |S_2|^2\Big\}\\[8pt]
&\lesssim T^{-2+C\ve}\Big\{ T^{-3\kappa+4}+\int_{[|t|>T^{\kappa/2}]} \min \Big(1, \frac Tt\Big)^2 |S_1|^2 |S_2|^2\Big\}
\end{aligned}
\eqno{(4.8)}
$$
where for the first term we used trivial bounds on $S_1$ and $S_2$.

It follows that the $\alpha$-parameter set where (4.7) fails is of measure at most
$$
T^{-3\kappa+C\ve}+ T^{-4+C\ve}  \int_{[|t|>T^{x/2}]} \min \Big(1, \frac T{|t|}\Big)^2 \big[ Av_{\alpha_1, \ldots, \alpha_{k-1}}|S_1(t)|^2] |S_2(t)|^2
dt\eqno{(4.9)}
$$
and it remains to bound the second term.

Estimate for $|t|>T^{\kappa/2}$,
$$
|S_2(t)|< T^{-\gamma} |I_k| \, |J_k|< T^{\frac 2k-\gamma}\eqno{(4.10)}
$$
for some $\gamma>0$ (a more precise bound will be unnecessary).
An estimate of the form (4.10) is obtained from standard exponential sum theory, exploiting only one of the variables $m$ or $n$.
Substituting (4.9), it remains to prove that

$$
\int_{[|t|>T^{\kappa/2}]} \min \Big(1, \frac T{|t|}\Big)^2 \big[Av_{\alpha_1, \ldots, \alpha_{k-1}} |S_1(t)|^2\big] dt\lesssim T^{4-\frac 4\kappa}\eqno{(4.11)}
$$
and (4.11) will clearly follow from

$$
\frac 1T\int\big[Av_{\alpha_1, \ldots, \alpha_{k-1}} |S_1(t)|^2\big] \vp \Big(\frac tT\Big) dt\lesssim T^{3-\frac 4k+\ve}\eqno{(4.12)}
$$
($0\leq \vp\leq 1$ a symmetric smooth bumpfunction).

The l.h.s. of (4.12) is bounded by
$$
\begin{aligned}
&\int_{[\frac 12, 1]^{k-1}}d\alpha_1\ldots d\alpha_{k-1}\\[10pt]
&\Bigg|\Bigg\{\begin{aligned} &(\bar m, \bar n, \bar {m'}, \bar {n'})\in (\mathbb Z^{k-1})^4; m_i, m_i' \in I_i, n_i, n_i' \in J_i \text { and}\\
&|n_1^k-m_1^k -(n_1')^k+(m_1')^k+\alpha_2 (n_2^k-m_2^k -(n_2')^k +(m_2')^k) +\cdots+ \alpha_{k-1} ( \ ) |< O(1)
\end{aligned} \Bigg\}\Bigg|.
\end{aligned}\eqno{(4.13)}
$$
Set $N=T^{\frac 1k}$ and
$$
N^\rho =\max(|n^k_1-m_1^k-(n_1')^k + (m_1')^k|,\ldots, | n^k_{k-1} - m^k_{k-1}- (n_{k-1}')^k + (m_{k-1}')^k|)
$$
$(0<\rho\leq k)$.
Then the (4.14) contribution to (4.13) is bounded by

$$
N^{-\rho}\big|\big\{ (x_1, x_2, x_3, x_4)\in \mathbb Z_+^4; x_i<N \text { and } |x_1^k- x_2^k+x_3^k- x_4^k|< N^\rho\big\}\big|^{k-1}.\eqno{(4.15)}
$$

To evaluate (4.15), we consider several cases for $\rho$.

If $\rho\geq k-1$ and since $k\geq 2$
$$
(4.15) \leq N^{-\rho} \Big(N^3 \frac{N^\rho}{N^{k-1}}\Big)^{k-1} \leq N^{3k-4}
$$
let $1\leq \rho< k-1$, $M\sim \max(x_1, x_2, x_3, x_4)=x_1$. Then

$$
\begin{aligned}
&\big|\big\{ (x_1, x_2, x_3, x_4)\in\mathbb Z_+^4; x_1\sim M, x_i<M \text { and } |x_1^k+x_2^k+x_3^k -x_4^k|< N^\rho\big\}\big|\leq\\
&M^3 \Big(1+\min\Big(M, \frac{N^\rho}{M^{k-1}}\Big)\Big) \leq N^3+N^{\frac {4\rho}k}+\frac {N^\rho}{M^{k-4}} 1_{[M^k>N^\rho]}
\end{aligned}
$$
and (4.15) contribution at most
$$
N^{3(k-1)-\rho}+N^{\frac {4\rho}k(k-1)-\rho} + N^{3k-4} \lesssim N^{3k-4}.
$$

Remains the case $\rho<1$.
If $k\geq 3$, estimate

$$
\begin{aligned}
&\big|\big\{ (x_1, \ldots, x_4)\in \mathbb Z^4_+; x_1\leq N \ \text { and } \ |x_1^k-x_2^k+x_3^k-x_4^k|<N^\rho\big\}\big|\leq\\[10pt]
&\int_{-1}^1 \Big|\sum_{x<N} e^{iux^k}\Big|^4 \ \Big|\sum_{n<N^\rho} e^{iun}\Big|du \ll N^{\frac 52+\frac \rho 2+\ve}
\end{aligned}
$$
using Cauchy-Schwarz and Hua's inequalities.
This gives the contribution
$$
N^{-\rho} N^{\frac {5+\rho}2(k-1)+\ve} \ll n^{3k-4+\ve}.
$$
For $k=2$, we obtain the contribution
$$
N^{-\rho} N^{2+\rho+\ve} \ll N^{2+\ve}= N^{3k-4+\ve}.
$$
This proves Theorem 2.

\end{document}